\documentclass[reqno]{amsart}

\pdfoutput = 1

\usepackage[english]{babel}
\usepackage{amsmath,amssymb}
\usepackage{amsfonts}
\usepackage{amsthm}
\usepackage{array,dcolumn}
\usepackage{graphicx}
\usepackage{algorithm}
\usepackage{algpseudocode}
\usepackage{todonotes}

\usepackage{tikz}
\usetikzlibrary{arrows,calc,decorations.pathmorphing}
\usepackage{pgfplots}

\usepackage[T1]{fontenc}
\usepackage[utf8]{inputenc}
\usepackage[activate={true,nocompatibility},spacing,kerning]{microtype}
\usepackage[hyperpageref]{backref}
\usepackage[pdftex,
colorlinks,
hyperfootnotes=false,
pdffitwindow=true,
plainpages=false,
pdfpagelabels=true,
pdfpagemode=UseOutlines,
pdfpagelayout=SinglePage,
pdftitle={An improved Talbot method for numerical Laplace transform inversion},
pdfauthor={Benedict Dingfelder and J.A.C. Weideman,
Technische Universität München and University of Stellenbosch},
pagebackref,
hyperindex,
filecolor=purple,
urlcolor=purple,
citecolor=blue,
linkcolor=blue]{hyperref}
\usepackage[sc,osf]{mathpazo}
\graphicspath{{plots/}}

\renewcommand{\leq}{\leqslant}
\renewcommand{\geq}{\geqslant}

\newtheorem{theorem}{Theorem}[section]

\begin{document}

\renewcommand{\figurename}{\small {\sc Figure\@}}
\renewcommand{\tablename}{\small {\sc Table\@}}

\title[An improved Talbot method for numerical Laplace transform inversion]{An improved Talbot method for numerical Laplace transform inversion}

\author{Benedict Dingfelder}
\address{\scriptsize{Zentrum Mathematik -- M3, Technische Universität München,
80290~München, Germany}}
\email{dingfelder@mytum.de}
\author{J.A.C. Weideman}
\address{\scriptsize{Department of Mathematical Sciences, University of Stellenbosch, Stellenbosch~7600, South
Africa}}
\email{weideman@sun.ac.za}

\begin{abstract}
The classical Talbot method for the computation of the inverse
Laplace transform is improved for the case where the transform
is analytic in the complex plane except for the negative real axis.   First, by
using a truncated Talbot contour rather than the classical
contour that goes to infinity in the left half-plane, faster convergence is achieved.
Second, a control mechanism for improving numerical stability is
introduced.
\end{abstract}

\maketitle

\section{Introduction}  \label{sec:intro}
Let $f(t)$ be real-valued and piecewise continuous
on $[0,\infty)$ and of exponential order $f = \mathcal{O}(e^{\gamma_0 t})$, $t \to \infty$,
for some real constant $\gamma_0$.  Then the Laplace transform of $f(t)$ is
defined by
\[
F(z) = \int_0^{\infty} e^{-zt} f(t) \, dt, \qquad \mbox{Re} \, z > \gamma_0.
\]
This paper deals with the numerical solution of the inverse problem,
i.e., given $F(z)$, compute $f(t)$ at a specified value of $t$.   The function
$F(z)$ is analytic in the half-plane $\mbox{Re}\,  z > \gamma_0$,
and if it can be evaluated in the complex plane (as opposed to just on the
real axis), this
analyticity can and should be exploited by good numerical methods.
One such method is numerical quadrature applied to the inverse formula
\begin{equation}
\label{eq:problem}
f(t) = \frac{1}{2\pi i} \int_{\gamma-i\infty}^{\gamma+i\infty} e^{zt} F(z) \, dz,
\qquad  \mbox{Re} \, \gamma > \gamma_0,
\end{equation}
known as the Bromwich integral.   Two effective quadrature rules are the simple trapezoidal
and midpoint rules, used in combination with contour deformation as first suggested
in \cite{But57}.

The purpose of the contour deformation is to exploit the exponential factor
in (\ref{eq:problem}).  In particular, assume the path of integration
in (\ref{eq:problem}) can be deformed to a Hankel contour, i.e., a contour
whose real part begins at negative infinity in the third quadrant, then winds around all singularities
and terminates with the real part again going to negative infinity in the second quadrant. Examples
are shown  in Figure~\ref{fig:Plot3Contours} below.   On such contours
the exponential factor causes a rapid decay, which makes the integral particularly
suitable for approximation by the trapezoidal or midpoint rules.
The contour deformation can be justified by Cauchy's theorem,
provided the contour remains in the domain of analyticity of $F(z)$.
Some mild restrictions on the decay of $F(z)$ in the left half-plane are also
required; for  sufficient conditions, see \cite{talbot,TWS06}.

Suppose such a Hankel contour can be parameterized by
\begin{equation}
C: \ z = z(\theta), \qquad -\pi \leq \theta \leq \pi,
\label{eq:contourC}
\end{equation}
where $\mbox{Re} \, z(\pm \pi) = - \infty$.   Then
\begin{equation}
\label{eq:deformed}
f(t) = \frac{1}{2\pi i} \int_{C} e^{zt} F(z) \, dz
= \frac{1}{2\pi i}  \int_{-\pi}^{\pi} e^{z(\theta) t} F(z(\theta)) z^{\prime}(\theta) \, d\theta.
\end{equation}
We approximate the latter integral by the $N$-panel midpoint rule with uniform
spacing $h = 2\pi/N$, which yields
\begin{equation}
\label{eq:midpointsum}
f(t) \approx \frac{1}{N i} \sum_{k = 1}^{N} e^{z(\theta_k) t} F(z(\theta_k)) z^{\prime}(\theta_k),
\qquad \theta_k = -\pi+\big(k-\frac12\big)h.
\end{equation}
If the contour (\ref{eq:contourC}) is symmetric with respect to the real axis,
and if $F(\overline{z}) = \overline{F(z)}$ (which is the case for $f(t)$ real-valued),
then half of the transform evaluations can be saved.  That is, one needs to consider
only quadrature nodes in the upper (or lower) half-plane.   When comparing
numerical results it is important to keep this point in mind as several other
papers, including \cite{duffy,talbotParam,WT07},  consider  quadrature
rules with a total of $2N$ nodes.

Popular contours $C$ are the parabola \cite{But57,GM05} and the hyperbola \cite{GM05,LP04,SST04},
as well as the cotangent contour  introduced by Talbot in \cite{talbot}.
Here we consider  the Talbot contour in the form
\begin{equation}
\label{eq:cotangentTalbot}
z(\theta) = \frac{N}{t} \, \zeta(\theta), \quad
\zeta(\theta) = - \sigma + \mu \, \theta \cot(\alpha \theta) + \nu \, i \, \theta,
\quad -\pi \leq \theta \leq \pi,
\end{equation}
where $\sigma$, $\mu$, $\nu$  and $\alpha$ are 
constants to be specified by the user.   In the
original Talbot contour the parameter $\alpha$ did not appear,
i.e., $\alpha = 1$.   For reasons outlined below, the
modification we propose is to consider $0 < \alpha < 1$, in
which case we refer to (\ref{eq:cotangentTalbot}) as the modified Talbot contour. The scaling
factor $N/t$ was introduced in the original paper \cite{talbot} and
also used in \cite{AV04,talbotParam}.   Assuming a fixed value of $t$,
this is the appropriate scaling for maximizing the convergence rate in
the trapezoidal/midpoint approximation as $N \to \infty$.


The practical issue addressed in this paper is the choice of
the parameters $\sigma, \mu, \nu$ and $\alpha$ in (\ref{eq:cotangentTalbot}).
In his original paper \cite{talbot}, Talbot suggested parameter
choices in the case $\alpha = 1$, for various configurations of the
singularities of the transform  $F(z)$.    On the other hand, in \cite{AV04}
a simpler approach was proposed that fixes
$\sigma = 0$, $\alpha = 1$, and $\mu = \nu = 0.2$ in  (\ref{eq:cotangentTalbot})
regardless of the singularity locations. (Note that the parameter $0.2$
is the constant $\frac25$ in the contour of \cite{AV04}, when
making allowance for the $N$ vs.~$2N$ convention as mentioned
below (\ref{eq:midpointsum})).

Existing software codes include  Algorithm 682, which is a Fortran package that implements
Talbot's original contour and parameter selection strategy
\cite{MR90}.   Alternatively, the parameter choices
of \cite{AV04} were incorporated in a Mathematica code \cite{abateValkoCode} as well as
a MATLAB code \cite{mcClureCode}.

These software implementations work well for many
transforms.   For the special case when all the
singularities of $F(z)$ are located on the negative real axis, however,
one can choose the parameters better, as argued in \cite{talbotParam}.
By using a saddle point analysis, optimal parameters for the
contour (\ref{eq:cotangentTalbot}), with $\alpha = 1$,  were derived that predict
a convergence rate of $\mathcal{O}(e^{-0.949 N})$
as $N \to \infty$ for  fixed $t >0$.  An application of the
same saddle point analysis to the contour
of  \cite{AV04} shows that the convergence rate in this case
is slightly weaker, namely $\mathcal{O}(e^{-0.676 N})$.  This
correlates well with the figure $\mathcal{O}(10^{-0.6M})$ obtained by
numerical experimentation in  \cite{AV04}, when $M = N/2$.

Somewhat surprisingly, simpler contours such as the parabola and
hyperbola actually yield faster convergence than Talbot's original
contour, namely $\mathcal{O}(e^{-1.047N})$
and $\mathcal{O}(e^{-1.176N})$ respectively; see \cite{WT07}.
However, if one introduces the parameter $\alpha$ in (\ref{eq:cotangentTalbot})
the convergence rate of the Talbot method
can be improved to $\mathcal{O}(e^{-1.358N})$, which is significantly
better than any of the other convergence rates cited.
The above convergence rates, converted to base-$10$, are summarized in
Table~\ref{tab:table}.

\begin{table}[htb] \label{tab:table}
\caption{Convergence rates cited in the literature. The estimates apply to transforms
whose singularities are restricted to the negative real axis, such as
$F(z) = (z+1)^{-1}$ or $z^{-1/2}$.
Here $M = N/2$, where $N$ is the total number of nodes in the quadrature
sum (\ref{eq:midpointsum}); in other words, $M$ is the number of transform
evaluations.   The hyperbola and the modified Talbot contour are the only
contours among these that gain at least one decimal digit per transform evaluation.}
\begin{center}
\begin{tabular}{ccc}
Contour      & Convergence  & Reference \\ \hline
  Talbot &  $10^{-0.6M}$  &  \cite{AV04} \\
  Talbot &  $10^{-0.8M}$  &  \cite{talbotParam} \\
  Parabola & $10^{-0.9M}$  & \cite{WT07}  \\
  Hyperbola & $10^{-1.0M}$ &  \cite{WT07}  \\
  Modified Talbot &  $10^{-1.2M}$ &  [this paper] \\
\end{tabular}
\end{center}
\end{table}

Optimal parameters that achieve this $\mathcal{O}(e^{-1.358N})$ convergence
rate were derived by one of the authors of the present
paper, and reported in \cite{TWS06}, later to be included
in text books such as  \cite[Sect.~5.3]{Strang}.  Details of the derivation
have not been published, however.   Having simplified and improved
the derivation recently, we now present the details here.

The main considerations are summarized in  Figure~\ref{fig:Plot3Contours}.
The best contour has to strike a balance between passing too close
to the singularities, or too far from them in which case the exponential
factor in (\ref{eq:deformed}) becomes too large.  The quadrature
nodes on the best contour should also extend just far enough into
the left half-plane to reach the desired accuracy but their contributions
should not be negligibly small.

\begin{figure}[htb]
\centering
\includegraphics[width=0.6\textwidth]{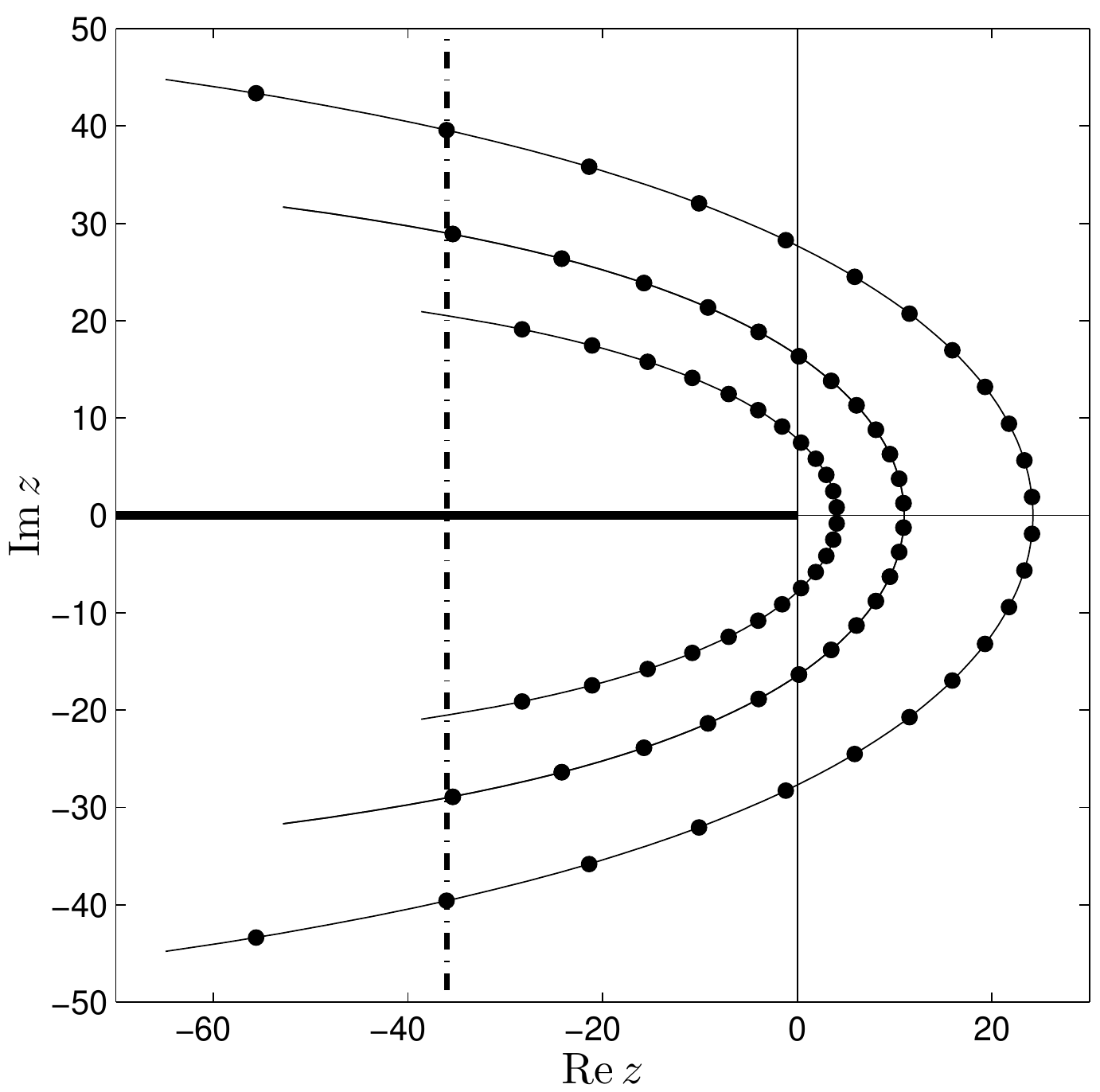}
\caption{Three possible Talbot contours (\ref{eq:cotangentTalbot}),
each shown with $N = 24$ nodes in the midpoint rule.   The
middle contour is close to optimal, in the case where $t=1$ and the
singularities of $F(z)$ are located on the negative real axis.
Its outermost conjugate pair of nodes reaches right up to the dash-dot curve,  where the
magnitude of the exponential factor in (\ref{eq:deformed}) becomes less than the
desired accuracy (here taken to be machine precision, $\epsilon \approx
2.2 \times 10^{-16}$).
The inner contour is suboptimal because (a) its outermost two nodes
do not reach quite up to the dash-dot curve and (b) it passes too
close to the singularities.   The outer contour is suboptimal because
(a) the contributions of the outermost pair of nodes are less than
the desired accuracy and therefore wasteful, and
(b) the exponential factor in (\ref{eq:deformed}) is too large in
the right half-plane.}
\label{fig:Plot3Contours}
\end{figure}

This discussion gives a clue as to why it is advantageous
to introduce the parameter $\alpha$ in (\ref{eq:cotangentTalbot}).
Because $\mbox{Re} \, z(\pm \pi) = -\infty$ in the original Talbot method,
too many nodes end up far out in the left half-plane where they make
no contribution to the quadrature sum (\ref{eq:midpointsum}).
We acknowledge  \cite{TrefPC} for this observation, and the
suggestion to consider a truncated Talbot contour instead.

On a technical note, an error is
committed when the truncated Talbot contour terminates at a finite location in the left half-plane.
In the original derivation \cite{talbot}, the path is closed at $\mbox{Re}(z)=-\infty$.
Because of the $e^{zt}$ factor, however, this error is of a similar magnitude as the truncation error---effectively
all contributions to the left of the dash-dot curve in Figure~\ref{fig:Plot3Contours}
are ignored.

The above  discussion also makes it clear that
the optimal contour parameters are determined by the
location of the singularities of $F(z)$.  Clearly
``one-size-fits-all'' parameter selections cannot be
optimal in all situations.   Below we continue to focus on the
case  when the singularities of $F(z)$ are located
on the negative real axis.  Although this restricts
the class of transforms significantly, this case is sufficiently
important in applications to justify such a study.
For example, when solving certain parabolic PDEs, singularities
on the negative real axis is an appropriate assumption
to make; see Section~\ref{sec:tests}.

%

We conclude by commenting on the scaling factor $N/t$ in (\ref{eq:cotangentTalbot}),
which has a few negative consequences.   First, because the nodes $z(\theta_k)$ depend
on $N$, function values cannot be re-used as $N$ is increased.
The increased convergence rate compensates for this, however.
Second,  for fixed $t$ the contour can
move far into the right half-plane as $N$ is increased.
Because of the exponential factor in (\ref{eq:problem}),
the terms in the summation can become large and hence
suffer from floating-point cancelation errors.
We address this issue in Section \ref{sec:roundoff}.   Third, because
the quadrature nodes depend on $t$ as well, the same
contour cannot be used for two or more distinct values
of $t$.  In this case the method would be inefficient when the transform
$F(z)$ is expensive to evaluate.   We shall not address
the latter issue here, but note that for the case $\alpha = 1$
it was considered in \cite{R95} and for the parabolic and hyperbolic
contours it was considered in \cite{WT07}.

The outline of the paper is as follows:  In Section~\ref{sec:error}
we present the main error analysis that leads to the proposed
values of $\sigma, \mu, \nu$ and $\alpha$.   The analysis is
similar to the one given in \cite{talbotParam}.  It is slightly
more intricate because of the additional parameter $\alpha$,
but at the same time we have made some improvements that
simplify the analysis.   A strategy for
improving numerical stability for large $N$ is discussed in
Section~\ref{sec:roundoff}.  Like the analysis of Section~\ref{sec:error},
it is restricted to the case when the singularities
of $F(z)$ are located on the negative real axis.
Section~\ref{sec:tests} contains the numerical tests, where
theory and experiment are compared.

 \section{Error estimates}
\label{sec:error}

Our basic error estimate is a well-known contour integral formula
for the error in the trapezoidal rule approximation of
functions that are analytic, periodic, and real-valued;
see for example  \cite[Thm.~9.28]{K}.  We retain the property
of analyticity but extend the theorem as follows:  (A) we consider
non-periodicity, (B) we consider complex-valued functions, and
(C) we consider the midpoint rule rather than the trapezoidal rule.
(A) is necessary because the integrand in (\ref{eq:deformed})  is not periodic.  (B) is necessary because
the integrand is complex-valued, and has very different
analytic properties in the upper and lower half-planes.
(C) is not really necessary, only a matter of convenience.
All three of these modifications can be readily incorporated into the proof
given in \cite[Thm.~9.28]{K}, which is an elementary application of the
residue theorem.  We omit the details.

\begin{theorem} \label{Kress} Let $\theta_k$ be defined as in (\ref{eq:midpointsum}).
 If $g: D \to \mathbb{C}$ is analytic in the rectangle
 $D = \{\theta \in \mathbb{C} : -\pi < \
\mbox{Re} \, \theta < \pi \mbox{ and }
 -d<\mbox{Im} \, \theta <c\}$ where $c,d>0$, then
\[ \int_{-\pi}^\pi g(\theta) \, d\theta - \frac{2\pi}{N} \sum_{k=1}^N g(\theta_k) =
\mbox{E}_+(\gamma) + \mbox{E}_-(\delta).
\]
Here
\[\mbox{E}_+(\gamma) = \frac{1}{2} \bigg(
\int_{-\pi}^{-\pi+i\gamma}  +
\int_{-\pi+i\gamma}^{\pi+i\gamma}  +
\int_{\pi+i\gamma}^{\pi}
\bigg)  \big(1+i \tan(\frac{N\theta}{2})\big)g(\theta) \, d\theta \]
and
\[\mbox{E}_-(\delta) = \frac{1}{2} \bigg(
\int_{-\pi}^{-\pi-i\delta}  +
\int_{-\pi-i\delta}^{\pi-i\delta}  +
\int_{\pi-i\delta}^{\pi}
\bigg)  \big(1-i \tan(\frac{N\theta}{2})\big)g(\theta) \, d\theta \]
for all $0<\gamma<c$, $0<\delta<d$ and $N$ even.  For odd $N$,
replace the $\tan(N\theta/2)$ with $-\cot(N\theta/2)$.
\end{theorem}

\begin{figure}[htb]
\begin{tabular}{cc}
\includegraphics[width=0.5\textwidth]{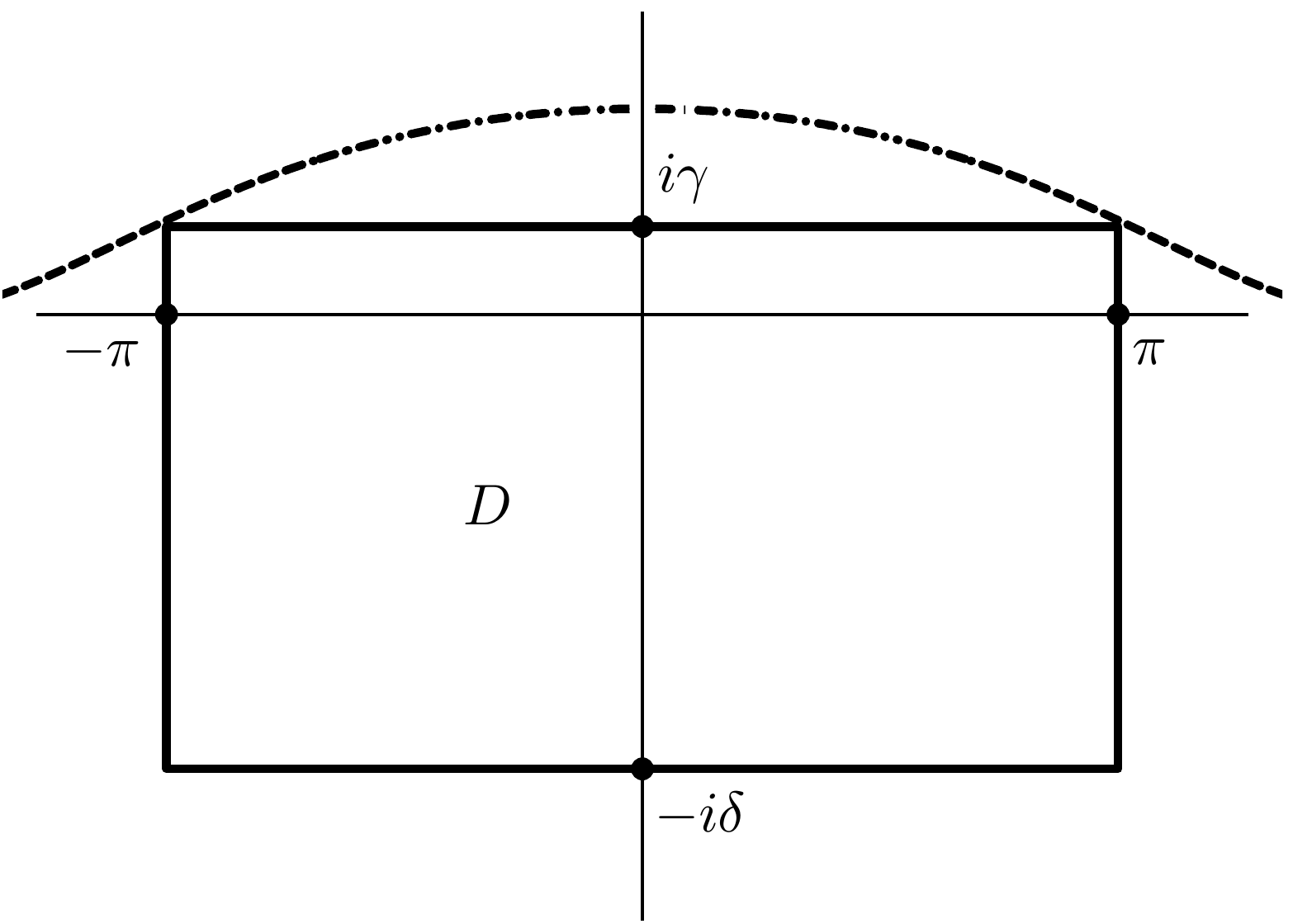} &
\includegraphics[width=0.5\textwidth]{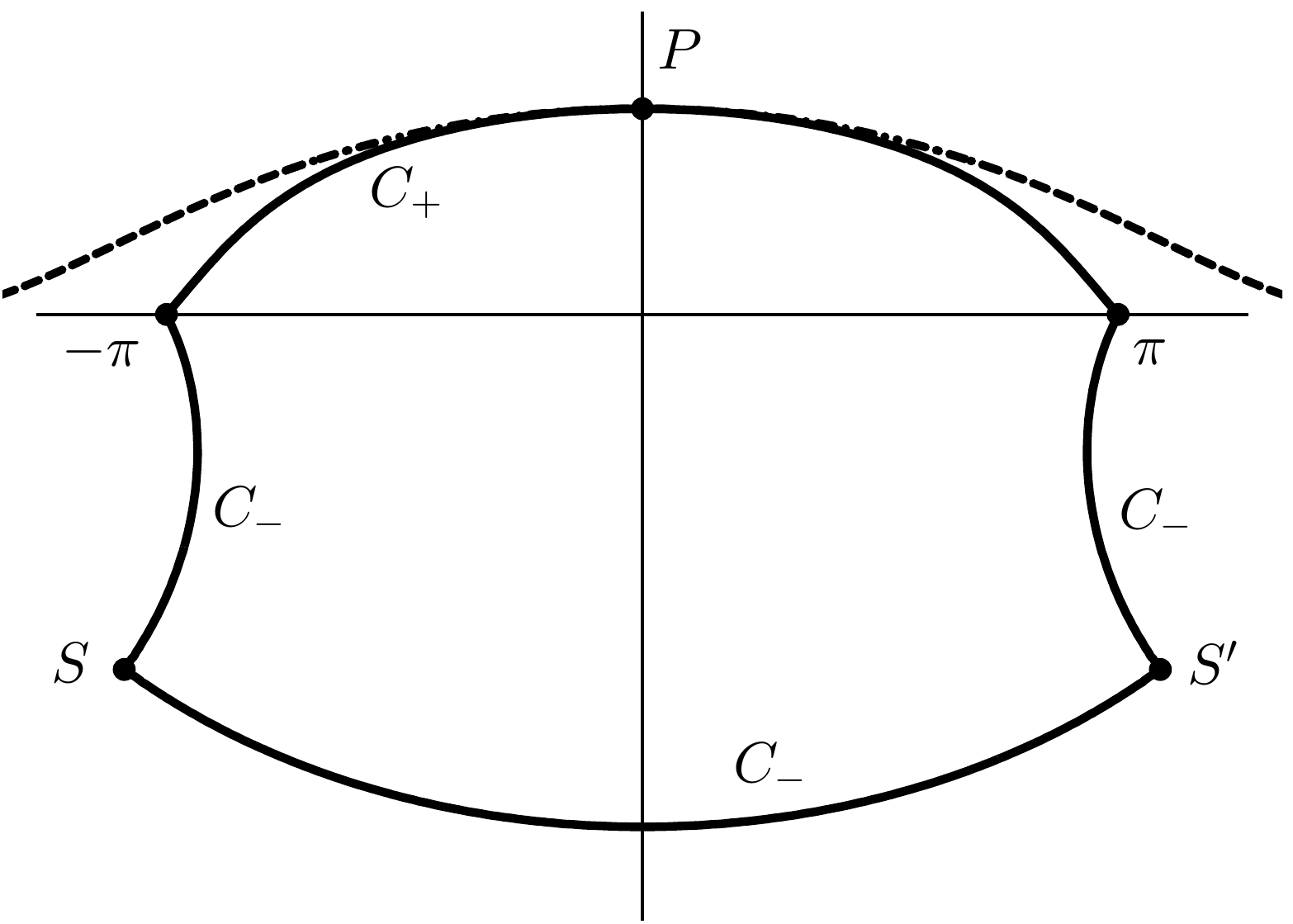}
\end{tabular}
\caption{Left:  Rectangular contour in the complex $\theta$-plane
for evaluating the error
integrals $\mbox{E}_+(\gamma)$ and $\mbox{E}_-(\delta)$
of Thm.~\ref{Kress}.  Right:  Deformed contour, along which
the magnitudes of the integrands of
$\mbox{E}_+(\gamma)$ and $\mbox{E}_-(\delta)$
are approximately constant.  The dash-dot curve in both diagrams
shows a possible image of the negative real axis of the $z$-plane,
where singularities might be located. The point $P$ is defined by the first two equations
in (\ref{eq:FirstThreeEquations}) and the points $S$ and $S^{\prime}$ by the first two equations in (\ref{eq:SecondThreeEquations}).
}
\label{fig:StripTalbot}
\end{figure}

The contours of integration for computing the error integrals
$\mbox{E}_+(\gamma)$ and $\mbox{E}_-(\delta)$
are shown in the left diagram of Figure~\ref{fig:StripTalbot}.
We remark that if $g(\theta)$ is real-valued on the real axis,
then $g(\overline{\theta}) = \overline{g(\theta)}$ and $d$ and $c$ can
be taken to be equal.  When it is moreover $2\pi$-periodic, the contributions
on the sides $\mbox{Re} \, \theta = \pm \pi$ cancel.  In that
case, if we define $\mbox{E}(\gamma) = \mbox{E}_+(\gamma) +
\mbox{E}_-(\gamma)$, then
\[
\mbox{E}(\gamma) = \mbox{Re}
\int_{-\pi+i\gamma}^{\pi+i\gamma}     \big(1+i \tan(\frac{N\theta}{2})\big)g(\theta)
\, d\theta.
\]
Here, and below, we assume $N$ even but the analysis for odd $N$ is similar.
By analyzing the behavior of the tangent function in the complex plane,
this can be bounded by \cite[Thm.~9.28]{K}
\[
| \mbox{E}(\gamma) | \leq \frac{4 \pi \mathcal{M}}{e^{cN}-1},
\]
where $\mathcal{M}$ is a bound on the magnitude of $g$ in the domain $D$.   For a given
value of $N$, two quantities
control the size of the error: $c$, which is the half-width
of the rectangle of analyticity, and $\mathcal{M}$, which is governed
by the growth of $g$ in the complex plane.

Here we apply the error formula of Theorem~\ref{Kress} to the case
\begin{equation}
g(\theta) = \frac{1}{2\pi i} \, e^{z(\theta)t} F(z(\theta)) z^{\prime}(\theta), \qquad -\pi \leq \theta \leq \pi,
\label{eq:gtheta}
\end{equation}
where $F(z)$ is the transform to be inverted.  This
is not a real-valued  function, which is why the rectangle
of integration  in Figure~\ref{fig:StripTalbot} is not shown
symmetric with respect to the real axis.  As we shall
see (and similar observations were made in \cite{talbotParam,WT07}),
in the upper half-plane the top edge of the rectangle is
restricted by the singularities of $F(z(\theta))$, while in
the lower half-plane the restriction is the size of
$e^{z(\theta)t}$.

To keep the analysis
tractable, we limit our discussion primarily  to the model transform
\begin{equation}
F(z) = (z+\lambda)^{-1}, \qquad \lambda > 0,
\label{eq:ModelTransform}
\end{equation}
i.e., singularities on the negative real axis only.   In
Section~\ref{sec:tests} we shall consider a matrix version
of this transform  when considering applications to PDEs.

For~(\ref{eq:ModelTransform}) and similar transforms, we can
estimate the error  using Theorem~\ref{Kress} and the rectangular
contour shown on the left in Figure~\ref{fig:StripTalbot}.
For a better estimate, one can try and deform the
contour as shown on the right in Figure~\ref{fig:StripTalbot}.

To explain how this can be done, consider the error integrals
in Theorem~\ref{Kress}.  Ignoring the constant factors $1/2$ in
those integrals, as well as the factor $1/(2\pi i)$
in (\ref{eq:gtheta}), the error formulas can be expressed as
\[
E_{\pm} = \int_{C_{\pm}} e^{N h_{\pm}(\theta)} \, d\theta,
\]
where
\[
h_{\pm}(\theta) = \zeta(\theta)  + N^{-1} \Big( \log\big(1\pm i \tan(\frac{N\theta}{2})\big)
+ \log F(z(\theta))  + \log z^{\prime}(\theta) \Big),
\]
and $C_+$ and $C_-$ refer to contours in the upper and lower half-planes, respectively.
Our strategy is now as follows: suppose it is possible to deform
the rectangle on the left in Figure~\ref{fig:StripTalbot}
to a contour along which $\mbox{Re} \, h_{\pm}(\theta)  = -c$ for some positive constant $c$.
Then
\begin{equation}
C_{\pm}: \ \mbox{Re} \, h_{\pm}(\theta) = -c \qquad \Longrightarrow \qquad |E_{\pm}| = O(e^{-cN}),
\label{eq:Cplusminus}
\end{equation}
and the value of $c$ can be maximized by picking the best among all such contours $C_{\pm}$.

The function $h_{\pm}(\theta)$ defined above is too complicated to analyze
in this manner,
and we make a few approximations in the case $N \gg 1$.
First, we drop the  terms involving $F(z(\theta))$ and
$z^{\prime}(\theta)$.   This means $F(z(\theta))$ should not be zero nor
exponentially small on any part of the contours $C_{\pm}$, and
likewise  $z^{\prime}(\theta)  \neq 0$.  Regarding the term involving
the tangent function in $h_{\pm}(\theta)$, we shall approximate it for
$N \gg 1$  as follows.  By converting the tangent function to exponentials
one obtains, with $\theta = x + iy$,
\[
 1 \pm i \tan\big(\frac{N}{2}(x+iy)\big)  = \frac{2 \, e^{\mp \frac12 N y} e^{\pm \frac12 i N x}}
{ e^{- \frac12 N y} e^{\frac12 i N x} + e^{\frac12 N y} e^{-\frac12 i N x}}.
\]
Taking the logarithm yields
\begin{equation}
\log\big(1\pm i \tan(\frac{N\theta}{2})\big) \sim  \mp N y, \qquad N \to \infty,
\qquad
\label{eq:est}
\end{equation}
valid for fixed, nonzero $y$.
The upper and lower sign choices correspond to $y>0$ and $y<0$, respectively.
Accordingly, we approximate the contour $\mbox{Re} \, h_{\pm}(\theta)  = -c$
in (\ref{eq:Cplusminus}) by
\begin{equation}
 C_{\pm}: \quad  \mbox{Re} \ \zeta(x+i y) \mp y = -c.
 \label{eq:ApproxCplusminus}
\end{equation}

First, consider the upper half of the $\theta$-plane.   As mentioned above, the
error contour $C_+$ is restricted by the location of the singularities
of $F(z(\theta))$ in the upper half of the $\theta$-plane.
For now, consider the model problem (\ref{eq:ModelTransform}).
Its pole is located at $z = -\lambda$, or
$\zeta = -\lambda t/N$.   In keeping with the limit $N \to \infty$,
$t$ fixed, we consider therefore $\zeta = 0$.  The same
argument will apply if there is a branch-cut on the negative
real $z$-axis, for in practice we need to consider only a finite section
of the branch-cut (e.g., the section between the origin and
the dashed line in Figure~\ref{fig:Plot3Contours}).

We therefore need to examine the zeros of $\zeta(\theta)$, i.e.,
\[
- \sigma + \mu \, \theta \cot(\alpha \theta) + \nu \, i \, \theta = 0.
\]
It is no easy task to establish theoretically how  the complex roots of this equation
depend on all possible parameter values.   We therefore resorted to numerical rootfinding, restricting
searches to the strip $-\pi \leq \mbox{Re} \, \theta \leq \pi$, which is   of interest here.
This revealed the following:  For certain parameter choices there are roots in the lower half
of the $\theta$-plane.
For other parameter choices the roots are in the upper half-plane, located
symmetrically  with respect to the imaginary $\theta$-axis.  For yet
other parameter values the roots are located precisely on the imaginary $\theta$-axis,
with the one closest to the real axis the critical one that limits
the contour $C_+$.  Using such numerical experimentation and keeping in mind that
$\zeta^\prime(\theta)$ may not be $0$ as remarked below (\ref{eq:Cplusminus}),
we came to the same conclusion as \cite{talbotParam}.  Namely,
the configuration that gives the widest domain of analyticity in the upper half-plane occurs
when the critical root (shown as the point marked $P$ in Figure~\ref{fig:StripTalbot})
is located on the positive imaginary $\theta$-axis and
is of double multiplicity.

This leads to the two  equations $\zeta(iy) = 0, \zeta^\prime(iy) =0$ for some $y>0$.
By considering (\ref{eq:ApproxCplusminus}), we conclude that $y = c$.
If we further let $\theta \to \pm \, \pi$
on the contour $C_+$, we get $\zeta(\pm \pi) = -c$, which is by symmetry just one equation,
not two.   In summary,    in the upper half-plane we require that
\begin{equation}
C_+: \quad \zeta(ic) = 0, \qquad \zeta^{\prime}(ic) = 0, \qquad
\zeta(\pi) = -c.
\label{eq:FirstThreeEquations}
\end{equation}

By using these three equations, it is possible to solve for
$(\sigma,\mu,\nu)$ in terms of   $(\alpha,c)$,
namely
\begin{equation}
\sigma = 2 \alpha c^2 B,\qquad
\mu = 2\sinh^2(\alpha c)  B, \qquad
\nu = (\sinh(2 \alpha c)-2\alpha c)B,
\label{eq:smn}
\end{equation}
where
\[
B = \frac{c \sin^2(\alpha\pi)}{2\alpha c^2 \sin^2(\alpha \pi)-\pi\sin(2\alpha\pi)\sinh^2(\alpha c)}.
\]

It remains to fix the parameters $\alpha$ and $c$.  For this,
we consider the lower half of the $\theta$-plane.   Here the
 contour $C_-$ is not restricted by analyticity but by the size
 of the exponential factor in (\ref{eq:deformed}).  Again numerical
 experimentation confirmed the findings of \cite{talbotParam},
 namely that there are two saddle points in the lower half of
 the $\theta$-plane that are critical.  (These are the points
 marked $S$ and $S^{\prime}$  in Figure~\ref{fig:StripTalbot}.)
 By making $C_-$ pass through these saddle points, we can
 enclose a region as big as possible.

 \begin{figure}[htb]
 \centering
 \includegraphics[width=0.95\textwidth]{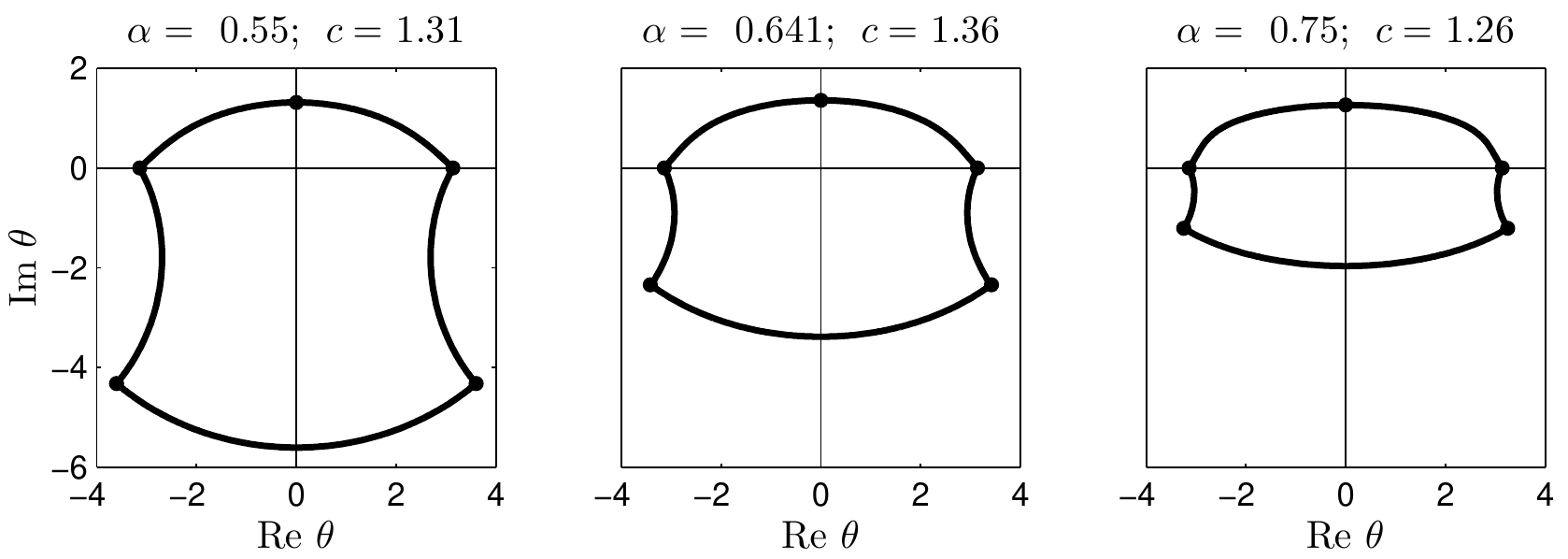}
 \caption{Three possible deformations of the rectangle of Figure~\ref{fig:StripTalbot}.
 The middle one maximizes the value of $c$, which is both the value of
the intersection on the positive imaginary axis and the decay rate in the error estimate
$E = \mathcal{O}(e^{-cN})$.}
\label{fig:StripTalbot3}
\end{figure}

 Denote the saddle point locations by   $\theta = x_s+i y_s$.  By
 inserting this into (\ref{eq:ApproxCplusminus}),  and taking partial derivatives, yields
 \begin{equation}
 C_-:  \quad
 \mbox{Re} \, \zeta^{\prime}(x_s+iy_s) = 0, \quad   \mbox{Im} \, \zeta^{\prime}(x_s+iy_s) = 1,
 \quad \mbox{Re} \ \zeta(x_s+i y_s) + y_s = -c.
 \label{eq:SecondThreeEquations}
 \end{equation}
We insert the
expressions for $\sigma$, $\mu$ and $\nu$ from (\ref{eq:smn}) and fix a value of $\alpha$.
Assuming a solution exists, the three equations (\ref{eq:SecondThreeEquations}) can then be solved
numerically for the unknowns $x_s$, $y_s$, and $c$.   By varying $\alpha$, the value of $c$
can thus be maximized.

We made no attempt to establish theoretical existence of a solution of the equations
(\ref{eq:SecondThreeEquations}), but direct numerical computation indicates that such a
solution exists for $\alpha$ in the interval  $[0.51, 0.82]$,
roughly.  Some of these configurations are shown in Figure~\ref{fig:StripTalbot3},  with the
middle one corresponding to the maximum value of $c$.    We computed this value  of
$c$ using a univariate optimization code, and found
\begin{equation}
\alpha = 0.6407; \qquad c = 1.3580.
\label{eq:AlphaC}
\end{equation}
Substitution of these two values into (\ref{eq:smn}) yields the modified Talbot contour
\begin{equation}
z(\theta) = \frac{N}{t} \big( -0.6122 + 0.5017 \, \theta \cot (0.6407 \, \theta) + 0.2645 \, i \, \theta \big),
\quad -\pi \leq \theta \leq \pi,
\label{eq:TalbotModified}
\end{equation}
as reported in \cite{TWS06}.   For completeness, we record that the corresponding
saddle points are located at $x_s+i y_s = \pm 3.4208 - 2.3438\, i$, correct
to all digits shown.

An alternative contour in which the cotangent is replaced by a rational function
was respectively suggested and analyzed in  \cite{talbot} and \cite{talbotParam}.
Applying the same strategy as outlined above, we derived
\begin{equation}
z(\theta) = \frac{N}{t} \big( 0.1446 +  {\frac{3.0232 \,  {\theta}^{2}}{{\theta}^{2}- 3.0767\,{\pi }^{2}}}+ 0.2339\,i \, \theta \big),
\quad -\pi \leq \theta \leq \pi.
\label{eq:TalbotModifiedRational}
\end{equation}
This rational contour is virtually as good as the cotangent contour, as it gives a convergence
rate of $\mathcal{O}(e^{-1.311 N})$, compared to the $\mathcal{O}(e^{-1.358 N})$ of
 (\ref{eq:TalbotModified}).  With either contour, an accuracy of close to machine precision
 $\epsilon \approx 2.2 \times 10^{-16}$ can be reached for values of $N$ around
 $26$ or $28$, i.e., about $13$ or $14$ transform evaluations.

We emphasize  that we have assumed in this section that all singularities of $F(z)$ are located on
the negative real axis, and $F(z)$ is not exponentially small as discussed below
(\ref{eq:Cplusminus}).      We also  ignored the effects of roundoff error, which
we consider next.

\section{Roundoff error control}
\label{sec:roundoff}

In the analysis of the previous section we explicitly assumed that $N \gg 1$, but tacitly
assumed infinite precision arithmetic.
When $N$ gets large, however, the contours (\ref{eq:TalbotModified}) and (\ref{eq:TalbotModifiedRational})
move far into the right half-plane.  The exponential factor in
(\ref{eq:midpointsum}) gets big, and there is loss of floating-point precision in
the summation.  Fortunately, for the theory of the previous section to be applicable,
$N$ needs to be only moderately large (we shall see that $N = 20$ or even
$N = 10$ is adequate).  There is therefore a range of values of $N$, say $N< N_*$,
where the error estimates of the previous section are accurate.   In the present
section we estimate the critical value of $N = N_*$ where roundoff error becomes
a factor, and we  propose an alternate set of parameters to improve stability
when $N > N_*$.

Using the same arguments as presented in \cite{Wei10} for the parabolic
contour, we model the error caused by finite precision by
\begin{equation}
R = \mathcal{O}\big(\epsilon \, e^{z(0)t} \big) = \mathcal{O}\big(\epsilon \, e^{N \zeta(0)} \big).
\label{RE}
\end{equation}
Here $\epsilon$ is the unit roundoff.

Let the implied constants in
(\ref{eq:Cplusminus}) and (\ref{RE}) be denoted by $k_1$ and $k_2$, respectively.
Then the  roundoff error becomes significant when
\begin{equation}
k_1  e^{-cN}  = k_2 \epsilon \, e^{N \zeta(0)}  \quad \Longrightarrow \quad
c + \zeta(0) + N^{-1} \log (\epsilon/k_0)  = 0,
\label{eq:Roundoff}
\end{equation}
where we defined $k_0 = k_1/k_2$.

Initially, let us suppose that the implied constants in (\ref{eq:Cplusminus}) and (\ref{RE})
are of the same order of magnitude, i.e., $k_0 \approx 1$.
Using the contour parameters as defined by (\ref{eq:smn})
and (\ref{eq:AlphaC}), a numerical solution of (\ref{eq:Roundoff}) with $\epsilon \approx  2.2 \times 10^{-16}$
yields that roundoff error becomes significant when $N \approx 24$ (25 in the case of the rational contour
(\ref{eq:TalbotModifiedRational})). An examination of the convergence curves shown in
Section~\ref{sec:tests} will confirm that $N = 24$ is often close to the critical value where
roundoff errors start to dominate. More generally, in a multiprecision environment with precision
equal to $d$ decimal digits, it follows from (\ref{eq:Roundoff}) that the critical value
of $N$ is approximately $1.5 d$.

Considering the general case where $k_0$ is not $O(1)$, this critical value of $N$ may be larger. We call this value
$N_*$, which can be estimated by analyzing the convergence behavior.  For $N \leq N_*$ roundoff error can be
neglected, i.e., the main contribution to the error
is due to the midpoint approximation. Hence the error behaves like $f(t)-f_N(t) = k_1  e^{-cN}$,
where $f_N(t)$ denotes the approximation (\ref{eq:midpointsum}). By computing $f_N(t)$ for a
range of values of $N$, it is straightforward to detect the value of $N$ at which the convergence
rate turns from exponential decay to exponential growth.
Having determined this critical value  $N = N_*$ we can insert it into
(\ref{eq:Roundoff}). By using the values of the contour parameters defined by (\ref{eq:smn})
and (\ref{eq:AlphaC}), it is possible to compute an approximate value of $k_0$.

We use this value of $k_0$ to compute new contour
parameters for any value of $N$, with  $N > N_*$.  Note that $\zeta(0) = -\sigma+\mu/\alpha$,
which depends on $c$ and $\alpha$ via (\ref{eq:smn}).
This means the equation  in (\ref{eq:Roundoff})
has two unknowns, so we fix $\alpha$ at the value (\ref{eq:AlphaC}).
Hence $c$ can be solved by applying a univaritate rootfinding
routine to the equation on the right in (\ref{eq:Roundoff}), for a given
value of $N (>N_*)$.  With $c$ and $\alpha$ known, the new contour parameters for
this particular value of $N$ are then computed
from (\ref{eq:smn}).

We can justify this approach by looking at the asymptotic behavior
of the parameters.   Letting $N \to \infty$ in (\ref{eq:Roundoff})
yields
\[
c = - \log (\epsilon/k_0) N^{-1} + O(N^{-3}), \qquad
\zeta(0) =  O(N^{-3}).
\]
Substitution of these two expressions into the error estimates
(\ref{eq:Cplusminus}) and (\ref{RE}), respectively, shows
that both reduce to $\mathcal{O}(\epsilon)$.   Hence errors
can be expected to remain at the roundoff level for
large $N$.

Because  $\zeta(0) = \mathcal{O}(N^{-3})$,
the contour is not only prevented from  moving too far into the
right half-plane,  it actually moves back to the left (and becomes
narrower) when $N > N_*$. This means that this strategy is only effective
for transforms $F(z)$ whose singularities are exclusively
located on the negative real axis. In the next section, we investigate
how well it works in practice.

We remark that a similar roundoff control strategy for the parabola was
considered in \cite{Wei10}. A related but different strategy for the hyperbola
was proposed in \cite{LPS}.  A new strategy that
automatically compute a contour of optimal stability for  (\ref{eq:problem})
was developed by one of the authors of the present paper.
Details can be found in \cite{dingfeld}.

\section{Applications and Comparisons}
\label{sec:tests}

In this section we verify the theoretical results of the previous sections by applying
the method to various model problems. These include a matrix transform that
can be used for the time integration of some parabolic PDEs,   as well as
two scalar transforms.  These two transforms were both taken from
applications and both were previously considered as test examples
in \cite{duffy}.   This reference also gives series or integral expressions
for the inverses, which we used  to compute  errors.

In all figures, the relative error is plotted against the number of nodes,
$N$,  in the midpoint  approximation (\ref{eq:midpointsum}).  For reference, we also plot the theoretically
predicted convergence rate $\mathcal{O}(e^{-1.358 N})$, as a dash-dot
curve.  We also remind the reader about the $N$ vs.~$2N$ convention mentioned
below (\ref{eq:midpointsum}).   When compared to the error curves
in  \cite{duffy,Wei10,WT07}, for example,
the values of $N$ on the horizontal axes of the graphs below should
therefore be halved.

We start with the model problem (\ref{eq:ModelTransform}) and its
matrix analogue
\begin{equation}
F_1(z; \lambda) = (z+\lambda)^{-1}, \qquad F_1(z; A) = (zI+A)^{-1}.
\label{eq:ModelTransform2}
\end{equation}
Here $\lambda>0$ and $A$ is a symmetric positive definite
matrix, with $I$ the identity matrix of the same size.

The matrix version  arises, for example, when
semi-discretizations of the heat equation are solved by Laplace
transform techniques \cite{GM05,LP04,SST04}.  Consider
\[
\mathbf{u}_t + A \, \mathbf{u} = \mathbf{0} \quad \Longrightarrow \quad
z \, \mathbf{U}(z) - \mathbf{u}_0  + A \, \mathbf{U}(z)  = \mathbf{0} \quad \Longrightarrow \quad
\mathbf{U}(z) = F_1(z; A)  \, \mathbf{u}_0,
\]
where $A$ is a $J \times J$ matrix representation of the negative Laplacian, $\mathbf{u}$(t) an $J \times 1$ vector of
unknowns with $\mathbf{U}(z)$ its Laplace transform, and $\mathbf{u}_0$ the initial condition.
The inverse formula (\ref{eq:problem}) yields
\begin{equation}
\mathbf{u}(t) = \exp(-At) \mathbf{u}_0 = \frac{1}{2\pi i} \int_{C} e^{zt}
 \big( F_1(z; A)  \, \mathbf{u}_0 \big) \, dz.
 \label{eq:MatrixExponential}
\end{equation}
In this case the singularities of the transform are the negatives of the eigenvalues $\lambda$
of $A$, which are located on the negative real axis because of the assumptions
on $A$.   The right-hand side is approximated
by the midpoint sum (\ref{eq:midpointsum}), where each quadrature node
requires the solution of a linear system with coefficient matrix
$\big(z(\theta_k)I + A\big)$ and right-hand side $\mathbf{u}_0$.

Our test example is  the heat equation $u_t - 0.01 \nabla^2 u = 0$ on $[0,1] \times [0,1]$,
supplemented with homogeneous Dirichlet boundary conditions.  We take
$A$ as the familiar block-tridiagonal matrix  based on the $5$-point
finite difference approximation to the negative of the Laplacian,
which is known to be positive definite.  The right-hand side $\mathbf{u}_0$
is taken to be random, and the reference solution
$\mathbf{u}(t)$ was computed using the spectral decomposition
of $A$, which is explicitly known.

\begin{figure}[htb]
\centerline{\includegraphics[width=0.9\textwidth]{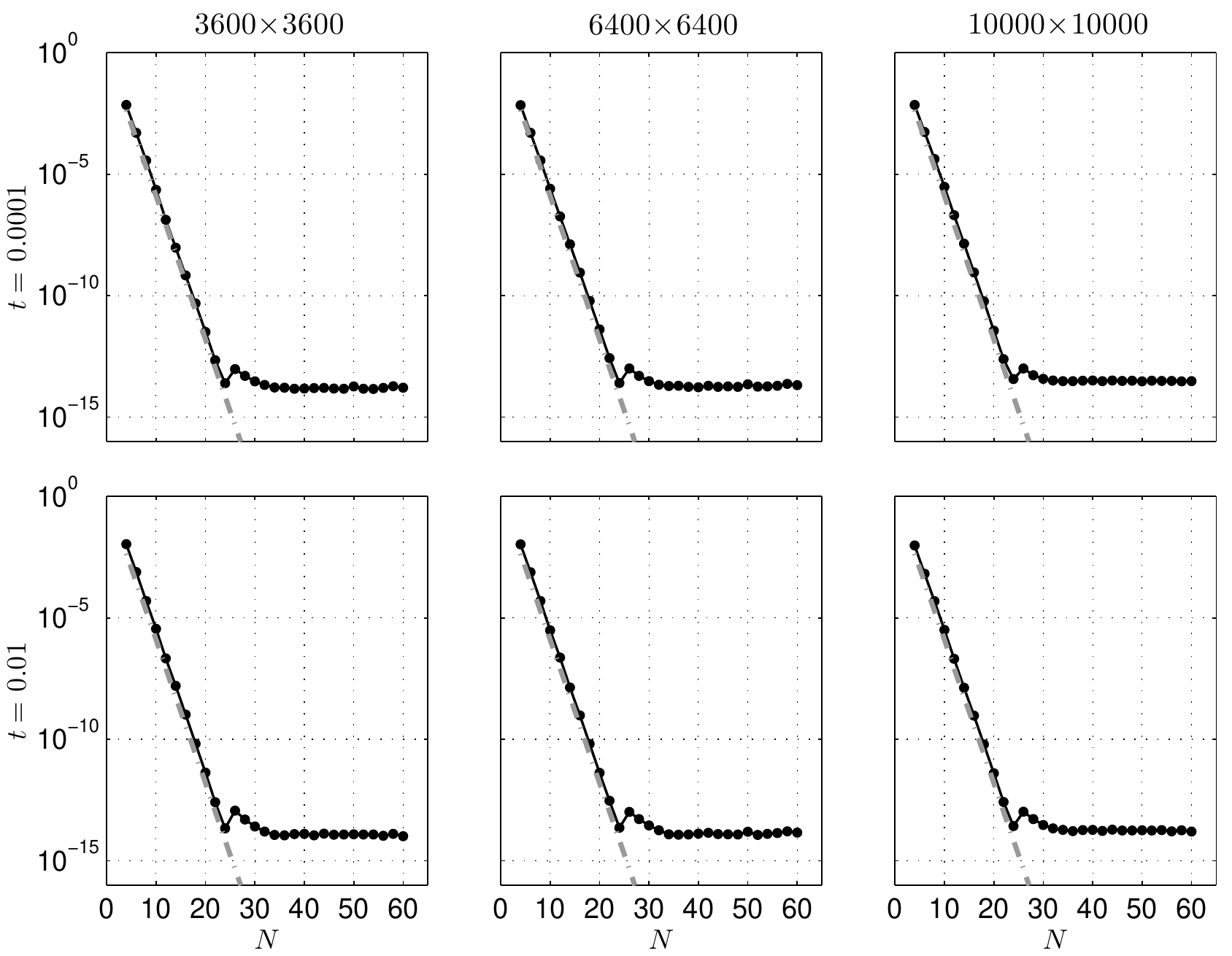}}
\caption{Relative errors in the $\infty$-norm when $\exp(-At) \mathbf{u}_0$
is approximated by applying the midpoint sum (\ref{eq:midpointsum}) to the
contour integral (\ref{eq:MatrixExponential}). For smaller values of $N$ the contour
(\ref{eq:TalbotModified}) was used, while for larger $N$ the roundoff
control procedure of Section~\ref{sec:roundoff} was used. The three columns
correspond to three sizes of $A$, while the two rows correspond
to two values of $t$. The dashed line represents the theoretical predicted convergence
rate $\mathcal{O}(e^{-1.358 N})$.}
\label{fig:Laplace5pt}
\end{figure}

The results for this example are shown in Figure~\ref{fig:Laplace5pt}.
Note that we have plotted the error here in the computed value
of $\exp(-At) \, \mathbf{u}_0$, i.e., only the temporal error.
The error in the actual solution of the PDE is not shown but it
would  not reach such small values  unless a
super accurate space discretization is used.

Figure~\ref{fig:Laplace5pt} confirms that the modified Talbot
contour converges as predicted by the error analysis of Section~\ref{sec:error}, namely
$\mathcal{O}(e^{-1.358 N})$ in all cases shown.    Ten-digit accuracy is reached with $N = 18$
or fewer nodes, i.e., no more than $9$ transform evaluations are required
(each of which involves the solution of a sparse linear system of the indicated dimension).
Almost full accuracy is reached for all $N > 24$.

The roundoff control strategy  described in the previous section is also seen to work well in practice.
To implement it, we set in the computations of Figure~\ref{fig:Laplace5pt} the value $k_0 = 1$ and
used the alternate parameters for all $N \geq 24$.

The remaining transforms in this section are scalar problems
from actual applications, as collected in   \cite{duffy}.
As an example of a transform with a series of poles on the
negative real axis we consider
\[F_2(z) = \frac{(100z-1)\sinh(\frac12 \sqrt{z})}{z(z\sinh(\sqrt{z})+\sqrt{z}\cosh(\sqrt{z}))}.\]
This is Test~$2$ in \cite{duffy}, where it was noted that the
apparent square root singularity can be removed by expanding the hyperbolic
functions.   We also consider Test~$3$ in \cite{duffy}, namely
\[F_3(z) = \frac{1}{z} \exp\bigg(-r\sqrt{\frac{z(1+z)}{1+cz}}\bigg),
\qquad r > 0, \ c > 0.  \]
This transform has a pole  at
$z=0$,  an essential singularity $z = -1/c$ as well as branch point
singularities at $z = 0$ and $z = -1$.   With  appropriate
definition of the branch cuts, this transform is analytic everywhere
off the negative real axis.

\begin{figure}[htb]
\centerline{\includegraphics[width=0.9\textwidth]{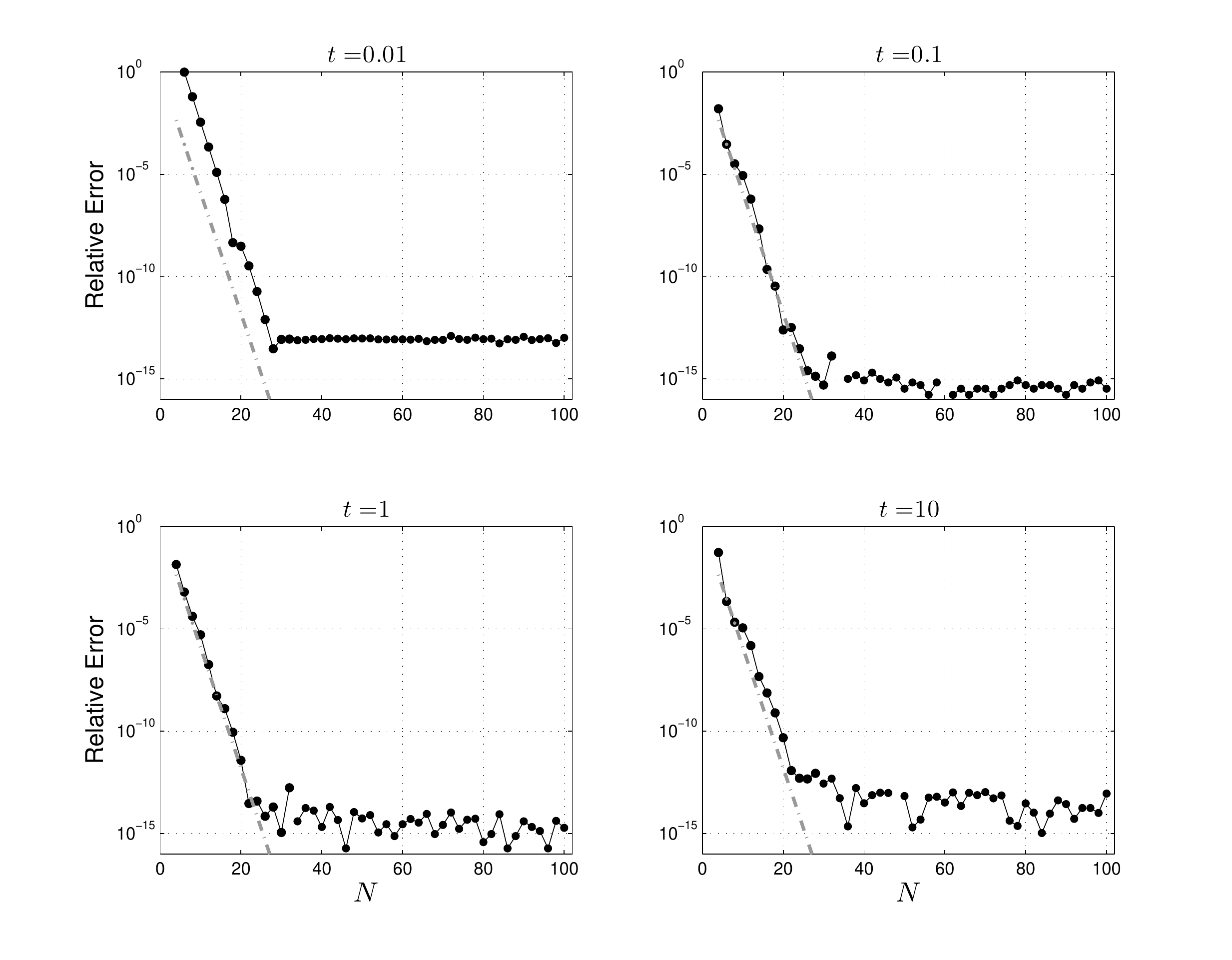}}
\caption{Relative error in the inversion of transform $F_2(z)$.
For smaller values of $N$ the contour
(\ref{eq:TalbotModified}) was used, while for larger $N$ the roundoff
control procedure of Section~\ref{sec:roundoff} was implemented.
}
\label{fig:Duffy2}
\end{figure}

\begin{figure}[htb]
\centerline{\includegraphics[width=0.9\textwidth]{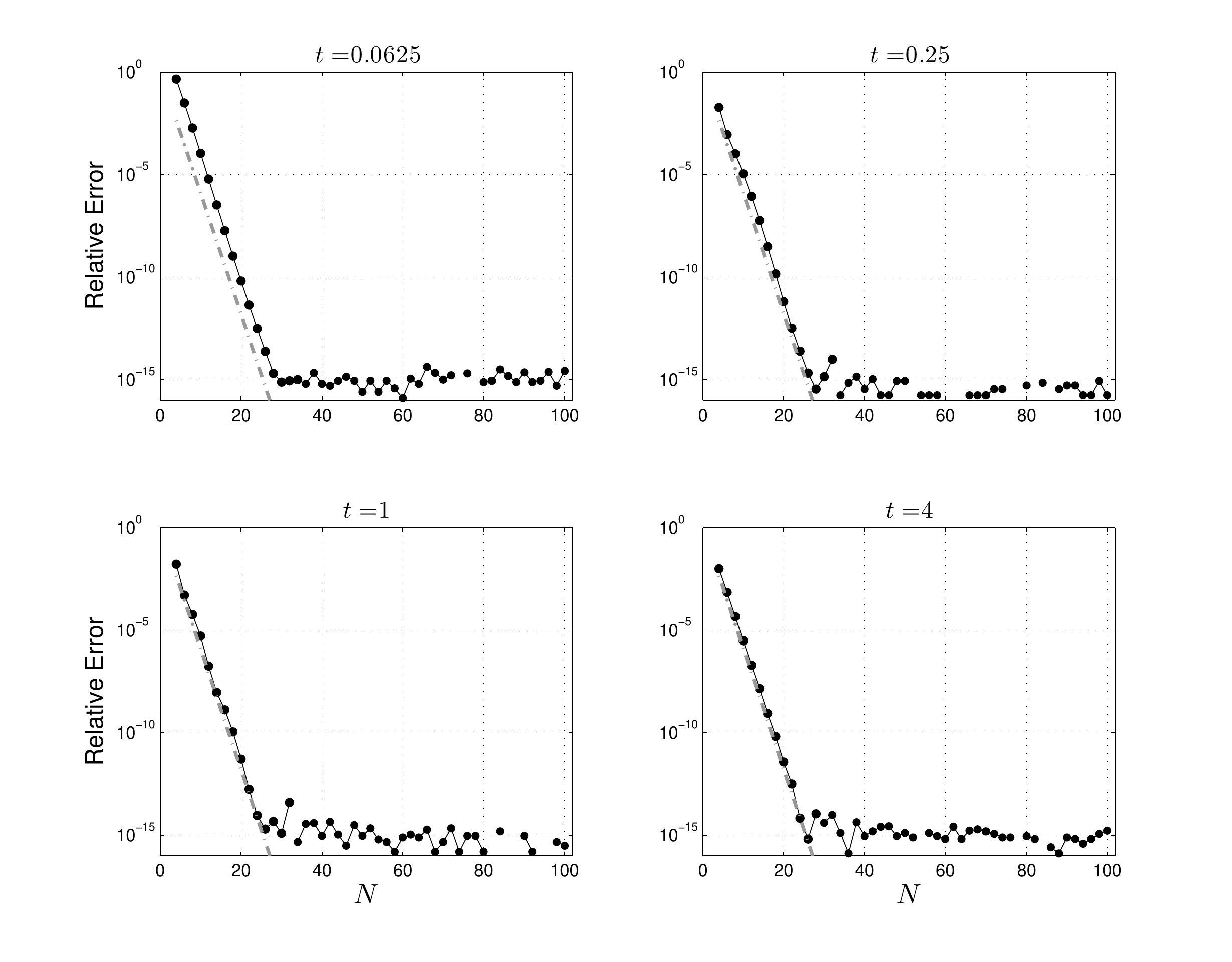}}
\caption{Same as Figure~\ref{fig:Duffy2} but the transform is $F_3(z)$,
with  $c = 0.4$ and $r = 0.5$.}
\label{fig:Duffy3}
\end{figure}

Numerical results for transforms $F_2(z)$ and $F_3(z)$ are
presented in Figures~\ref{fig:Duffy2} and \ref{fig:Duffy3}, respectively.
(The values of $t$ in these figures were chosen to match those in \cite{duffy}.)
Again we see empirical confirmation of the theoretical $\mathcal{O}(e^{-1.358 N})$ convergence
rate and the effectiveness of the roundoff control strategy.

A comparison with the results presented in \cite{duffy} shows a faster rate of convergence here.
For example, to achieve ten-digit accuracy at $t=1$, the midpoint sum (4) requires
$N = 18$, i.e., 9 transform evaluations, for both transforms  $F_2$ and $F_3$.
In contrast, the algorithm used in \cite{duffy}  achieves ten-digit accuracy
using 17 and 13 transform evaluations for $F_2$ and $F_3$, respectively.
The corresponding number of transform evaluations for the method of \cite{AV04}
is 15 and 14.

\begin{figure}[htb]
\centerline{\includegraphics[width=0.9\textwidth]{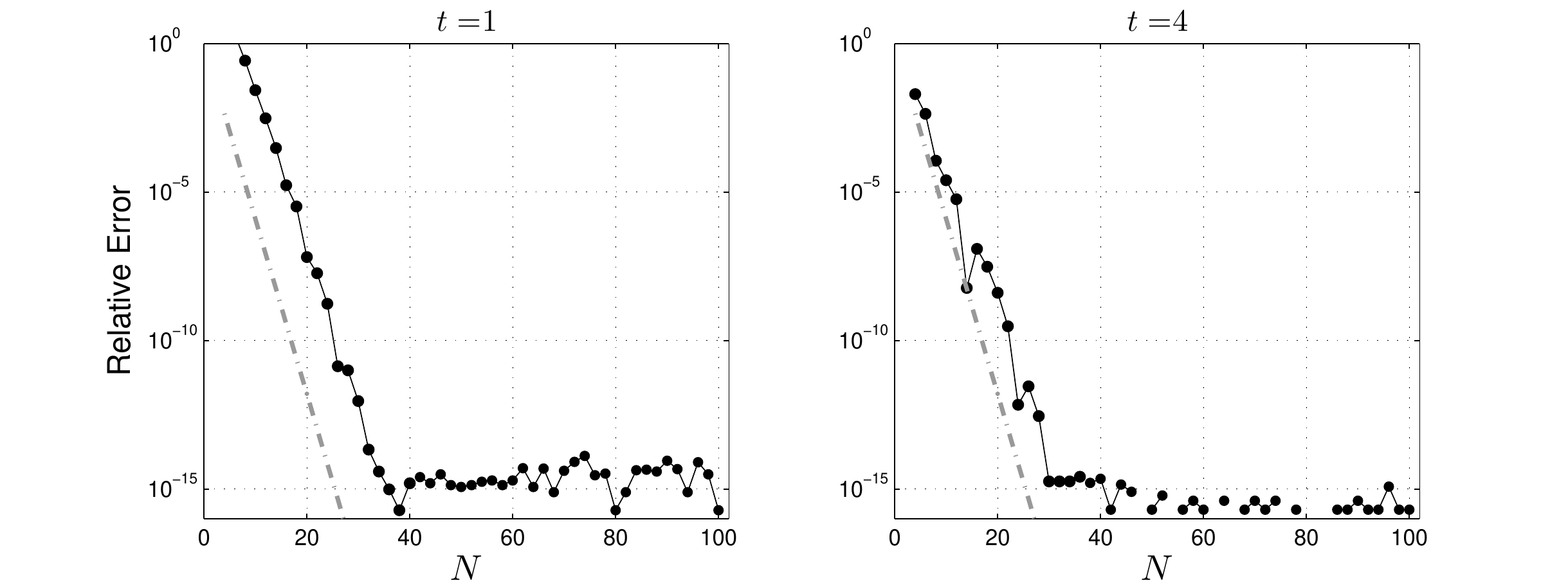}}
\caption{Same as Figure~\ref{fig:Duffy2} but the transform is $F_3(z)$,
with  $c = 0.4$ and $r = 3$.}
\label{fig:Duffy3a}
\end{figure}

Note that $F_3(z)$ is an example for which the error analysis of
Section~\ref{sec:error} can be less reliable.  The reason is its rapid
decay as $z \to +\infty$, particularly if $r \gg 1$; recall the discussion below (\ref{eq:Cplusminus}).
To see what happens, we show in Figure~\ref{fig:Duffy3a} convergence curves in the case
$r = 3$.    Asymptotically, the convergence rate is still  $\mathcal{O}(e^{-1.358N})$
(as indicated by the dash-dot curve)  but the implied constant
has increased by some orders of magnitude.     This also means that the critical value
$N$ where roundoff error becomes significant has now increased from $N = 24$
to about $N = 38$ and $N = 32$, respectively, in the cases $t = 1$ and $t=4$.
Nevertheless, the roundoff control strategy of Section~\ref{sec:roundoff}
is reasonably successful in preventing any catastrophic error growth.

\section{Conclusions}
\label{sec:conclusions}

A truncated version of Talbot's famous contour for numerical Laplace
transform inversion was analyzed.   For the case where the singularities
of the transform are located on the negative real axis, a new type of
error analysis was proposed and used to derive optimal contour
parameters.   With these parameters an exponentially convergent
method is obtained, with decay constant better than other known
contours.    In addition, a roundoff control strategy was proposed
for improved numerical stability when the number of nodes in
the quadrature scheme gets  large.    Numerical experiments
on transforms taken from actual applications confirmed the accuracy of
the optimal error estimates as well as the efficiency of the
roundoff control scheme.

We believe the improvements to Talbot's method suggested in this paper should
be considered for inclusion in future versions of software such as those
mentioned in Section~\ref{sec:intro}.
Our own MATLAB code that was used in generating
the numerical results of Section~\ref{sec:tests} can be downloaded
from \cite{Soft}.

\bigskip

{\bf Acknowledgment.}
The research of BD was supported by the DFG Collaborative Research Center TRR 109,
``Discretization in Geometry and Dynamics.''  He further acknowledges support
from the graduate program TopMath of the Elite Network of Bavaria and the
TopMath Graduate Center of TUM Graduate School at Technische Universit\"at M\"unchen.
The research of JACW was supported the National Research Foundation of South Africa.

Folkmar Bornemann and Nick Trefethen both offered valuable suggestions,
as have two anonymous referees.

\end{document}